\documentclass[11pt]{article}

\usepackage[english,russian]{babel}
\usepackage[T2A]{fontenc}
\usepackage[cp1251]{inputenc}
\usepackage{amssymb}
\usepackage{amsmath}

\usepackage{amsfonts,eucal,bm,amscd,amssymb,amsmath,latexsym}

\newcommand {\im} {\mathop {\rm im} \nolimits}

\begin{document}

\def\bibname{References}
\begin{center}
{\small\bf  ON VERSION OF FARKAS LEMMA
OF ALTERNATIVE
} \\ [2mm]
A.I. Golikov\\[2mm]
{\small\it Dorodnicyn Computing Centre of RAS}
\end{center}

\begin{center}
{\small\parbox[t]{10cm}{\hspace*{7mm}
The paper provides a version of Farkas lemma of alternative linear systems, when the alternative systems having different matrices of various number of dimentions.\\
\hspace*{7mm}{\bf Key words: }  Farkas lemma, alternative linear systems, unconstrained optimization}}
\end{center}

Let the system determining a set $X$ take the form of
$$
     Ax= b,~~x\ge 0_n.
\eqno (I)
$$
where $A$ is matrix $m\times n$, vector $b \in R^m$, $\|b\|\ne 0$.

The alternative system determining set $U$ can be presented as
$$
 A^\top u\le 0_n,~~b^\top u= \rho > 0,
\eqno (II)
$$
where $\rho $ --- an arbitrary fixed positive constant.

One, and only one of these systems, either ($I$) or ($II$), is always consistent, but never both.
 In case expression $b^\top u>0$ is used for ($II$) the above statement is known as Farkas lemma \cite{1}.

Besides their theoretical importance the theorems of alternative are of considerable value for computations \cite{2}.

For a given linear system, an alternative system is constructed in the space whose dimension
is equal to the number of equations and inequalities in the original system (not counting
constraints on the signs of variables). The original solvable system is solved by minimizing
the residuals of the inconsistent alternative system. The results of this minimization are used
to find the normal solution (with a minimal Euclidean norm) to the original system.

The replacement of the original problem by the minimization of the residuals of the inconsistent
alternative system may be advantageous when the dimension of the new variables is less than
that of the starting ones. In this case, such a reduction results in the minimization problem in a
space of lower dimension and allows one to obtain the normal solution to the original problem \cite{2}.

One and the same matrix $A$ and vector $b$ have always been used in alternative linear systems.
The paper shows a different way of alternative systems involving application of different matrices with various dimensions, which can be advantageous from computational point of view.

To determine a system resolvability and to find a solvable problem solution it suffices to find just a single vector $x^*$ or $u^*$ of the below-mentioned problems of quadratic minimization over the positive orthant or unconstrained minimization of a piecewise quadratic function
\begin{eqnarray}
\label{7.1}
&\min \limits_{x\in R^n_+}\frac 1 2\|b-Ax\|^2=\frac 1 2\| b- Ax^*\|^2,\\
\label{7.2}
&\min \limits_{u\in R^m}\frac 1 2\{\|(A^\top u)_+\|^2+(\rho-b^\top u)^2\}= \\
&=\frac 1 2\{\|(A^\top u^*)_+\|^2+(\rho-b^\top u^*)^2\}.\nonumber
\end{eqnarray}

The problems mutually dual  to  ($\ref{7.1}$) and ($\ref{7.2}$) respectively  will be the below-listed strictly concave quadratic programming  problems
\begin{eqnarray}
   \label{2.D1}
   \max \limits_{z \in Z}\{ b^\top z-\frac 1 2 \|z\|^2\},~~Z=\{z\in R^m~:~
   A^\top z\le 0_n\},
\end{eqnarray}
\begin{eqnarray}
   \label{2.D2}
&   \max \limits_{w\in W}\{\rho w_2-\frac 1 2\|w_1\|^2 -\frac 1 2w_2^2\},\\
& W=\{w_1\in R^n_+, w_2
   \in R^1~:~Aw_1-bw_2=0_m\} \nonumber.
\end{eqnarray}

From the dual features it follows that the solution $z^*$ to problem (\ref{2.D1}) can be expressed through problem (\ref{7.1}) solution by equation $z^*=b-Ax^*$; taking this equation into account, one can obtain $\|z^*\|^2=b^\top z^* $ from the equality of the objective functions' optimal values.

Similarly the solution $w_1^*$, $w_2^* $ to problem (\ref{2.D2}) can be expressed through the solution of problem (\ref{7.2}) in the following way: $w_1^*=(A^\top u^*)_+$, $w_2^*=\rho -b^\top u^*$ and $\|w_1^*\|^2+{w^*_2}^2=\rho w_2^*$ takes place.

If $X\ne \emptyset$, then $w_2^*> 0$ and normal (with minimal Euclidean norm) solution to system ($I$) can be expressed through the solution of (\ref{7.2}) as follows:
\begin{eqnarray}
   \label{7.4}
\tilde x^*=(A^\top u^*)_+/(\rho-b^\top u^*)=w_1^*/w_2^*.
\end{eqnarray}

If $X=\emptyset$, then $\|z^*\|\ne 0_m$, the normal solution to system ($II$) takes the form
\begin{eqnarray}
   \label{7.6}
   \tilde u^*=\rho (b-Ax^*)/\|b-Ax^*\|^2.
\end{eqnarray}

  The consideration below represents a special case of system ($I$) where matrix $A$ has rank $m$, i.e. $m \le n$. For the case concerned it will be shown that the system alternative to ($I$) can take a form different from ($II$), i.e. the alternative system can incorporate a matrix differing from $A$ and a vector other than $b$.

If $m\le n$ then the system
\begin{eqnarray}
   \label{7.7}
             Ax=b
\end{eqnarray}
is always solvable but its solutions may fail to include any nonnegative ones.
Let  $\bar X$ denote the set of system (\ref{7.7}) solutions.
Note that set $\bar X$ is always nonempty in contrast to set $X$.
The general solution of the  system of linear  equations (\ref{7.7}) can be
written in the form
\begin{eqnarray}\label{OR}
x= \bar x -K^\top y,
\end{eqnarray}
where $\bar x$ is a particular  solution of the system, and $K^\top y$ is the general solution
of the homogeneous system $Ax=0_m$, and $y\in R^\nu$. The  matrix $K$ can be chosen to be any $(\nu \times n)$
matrix such that its $\nu $ rows form a basis of the null space of
$A$ where $\nu=n-m$ is the defect of matrix $A$.
Therefore, $AK^\top=0_{m\nu}$. Here $0_{ij}$ denotes $(i\times j)$ matrix with zero entries.

Matrix $K$ is not uniquely defined; it can be constructed in various ways. If we partition
the matrix $A$ as $A=[B\,|\,N]$, where $B$ is non\-degenerate, then we can represent $K$ as $K=[-N^\top(B^{-1})^\top|\,I_\nu]$.
If we reduce $A$ by means of Gauss--Jordan transformations to the form $A=[I_m\,|
\,N]$, then we can represent $K$ as $K =[-N^\top|\,I_\nu]$ \cite{4}.

Let us determine the set $Y$ as
\begin{eqnarray}\label{Y}
Y=\{y\in R^\nu ~:~ \bar x - K^\top y \ge 0_n\}.
\end{eqnarray}

Equation (\ref{OR}) can be considered as an affine mapping
from $R^\nu $ to $R^n $. Here the image of set $Y$ is set $X$
specified by system ($I$). There exists a one-to-one correspondence between $X$ and $Y$.

     Indeed, for any $y \in Y$ equation (\ref{OR}) uniquely determines $x\in X$, i.e.
\begin{eqnarray}\label{X=Y}
X=\bar x-K^\top Y
\end{eqnarray}
In case of a full-range overdetermined  system (\ref{OR}) containing $n$ linear equations and
$\nu $ variables $y$ a pseudosolution
\begin{eqnarray}\label{!!}
y=(KK^\top )^{-1}K(\bar x - x)=(K^\top )^+(\bar x -x),
\end{eqnarray}
always exists. It solves (\ref{OR}) and is unique if and only if $\bar x - x \in \im K^\top $. This inclusion holds
if and only if $x \in \bar X$. Thus, for any $x \in \bar X$, formula (\ref{!!}) determines an affine transformation
that is the inverse of (\ref{OR}). Therefore, one can write
\begin{eqnarray}\label{YY}
Y=(K^\top )^+(\bar x -X).
\end{eqnarray}

So the following two systems
$$
  Ax=b,~~~x\ge 0_n, \eqno (I)
$$
$$
K^\top y \le \bar x, \eqno (I_y)
$$
are either  simultaneously solvable and interconnected by expressions (\ref{X=Y}) и (\ref{YY}) or simultaneously unsolvable if there exist no nonnegative general solution  $x= \bar x -K^\top y$ to system ($I$).

By  Gale theorem \cite{1} the following system determining set $V$ will be alternative
to system ($I_y$),
$$
Kv=0_\nu,~~~-\bar x^\top v=\rho >0,~~~v\ge 0_n. \eqno (II_v)
$$
System ($I_y$) being equivalent to system ($I$), system ($II_v$)
is simultaneously alternative to system ($I$).

The general solution to homogeneous system $Kv=0_\nu$ can be expressed by
matrix $A$ as $v=-A^\top u$. By changing the variables $v=-A^\top u$,
one can present system ($II_v$) as follows:
$$
A^\top u\le 0_n,~~~b^\top u=\rho >0.
\eqno (II)
$$
System ($II$) is alternative to ($I$), hence to ($I_y$).

If set $V$ is nonempty then set $U$ determined by system ($II$) is nonempty too, the two sets
having a one-one mapping expressed by:
$$
V=-A^\top U,~~~U=-(A^\top )^+V,
$$
where pseudoinverse matrix $(A^\top )^+$ is as follows:
$(A^\top )^+=(AA^\top )^{-1}A$.

The alternative systems interrelation can be represented as follows:

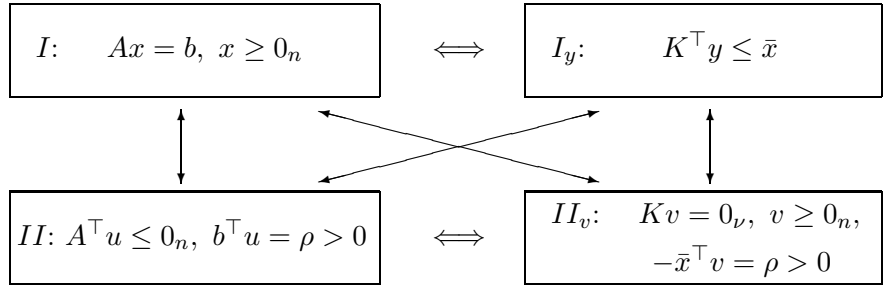
\begin{center}
\begin{picture}(350,35)
\put(0,0){\line(1,0){140}}
\put(0,35){\line(1,0){140}}
\put(0,0){\line(0,1){35}}
\put(140,0){\line(0,1){35}}
\put(10,15){$I$:~~~~~$Ax=b,~x\ge 0_n$}
\put(160,15){$\Longleftrightarrow$}
\put(195,0){\line(1,0){135}}
\put(195,35){\line(1,0){135}}
\put(195,0){\line(0,1){35}}
\put(330,0){\line(0,1){35}}
\put(205,15){$I_y$:~~~~~~~~$K^\top y\le \bar x$}
\end{picture}

\begin{picture}(350,70)
\put(0,0){\line(1,0){140}}
\put(0,35){\line(1,0){140}}
\put(0,0){\line(0,1){35}}
\put(140,0){\line(0,1){35}}
\put(3,15){$II$:~$A^\top u \le 0_n,~b^\top u=\rho >0$}
\put(160,15){$\Longleftrightarrow$}
\put(65,52){\vector(0,1){14}}
\put(65,52){\vector(0,-1){14}}
\put(195,0){\line(1,0){135}}
\put(195,35){\line(1,0){135}}
\put(195,0){\line(0,1){35}}
\put(330,0){\line(0,1){35}}
\put(205,23){$II_v$:~~~~$Kv=0_\nu ,~v\ge 0_n,$}
\put(243,5){$-\bar x^\top v=\rho >0$}
\put(265,52){\vector(0,1){14}}
\put(265,52){\vector(0,-1){14}}
\put(170,52){\vector(4,1){53}}
\put(170,52){\vector(-4,1){53}}
\put(170,52){\vector(4,-1){53}}
\put(170,52){\vector(-4,-1){53}}
\end{picture}
\end{center}
\bigskip

The double arrows correspond to simultaneously solvable/unsolvable systems
 and the ordinary ones stand for alternative systems.

Let us provide linear programming interpretation of Farkas lemma.
Here system ($I$) can be presented as the primal linear programming  problem
with its objective function coefficient vector identically equal to zero.
$$
\min_{x\in R^n_+}~\{0^\top_n x~:~Ax=b,~x\ge 0_n\}.
\eqno (P)
$$
The problem dual to ($P$) is as follows:
$$
\max_{u \in R^m}~\{b^\top u~:~A^\top u \le 0_n\}. \eqno (D)
$$

     It is common knowledge that for any couple of primal and dual LP problems there always
exists one of the following four cases:

     1) both primal and dual problems have  solutions;

     2) a primal problem is inconsistent and the dual one is unboarded;

     3) a primal problem is unboarded  and the dual one is inconsistent;

     4) both primal and dual problems are inconsistent.

For problems($P$) and ($D$) the latter two conditions cannot be fulfilled
because the constraints in ($D$) are always consistent, vector $u=0_m$ is feasible.

The first two cases are only possible.

In case 1) the optimal values of goal functions for problems ($P$) and ($D$)
are equal to zero and inequation $b^\top u\le 0$ holds for all feasible vectors $u$
owing to the weak duality theorem.
 Hence it follows solvability of system ($I$) and unsolvability of system $(II)$
$$
A^\top u \le 0_n, ~~b^\top u=\rho >0. \eqno (II)
$$

 In case 2) system ($I$) is inconsistent and system ($II$) is consistent for any $\rho >0$ due to
 unboundedness of dual problem ($D$).

 So one can obtain the simplest proof that  ($I$) and ($II$) are  alternative system
  employing a specific type of linear programming problems ($P$) и ($D$) and
  linear programming duality theory.

Problem ($\ref{7.1}$) can be considered an auxiliary problem of penalty function method as applied to
problem ($P$). Problem ($\ref{7.2}$) can be treated an auxiliary problem of  Morrison method with its parameter being equal to $\rho$
when applied to problem ($D$).
\bigskip

Acknowledgment. The work was supported by the Russian Foundation for Basic Research, project 15-01-08259.

\end{document}